\documentclass[12pt]{article}
\usepackage[top=2.54cm, bottom=2.54cm, left=2.70cm, right=2.70cm]{geometry}
\usepackage{tikz}
\usetikzlibrary{arrows.meta, decorations.markings, calc}
\usepackage[tbtags]{amsmath}
\usepackage{amssymb}
\usepackage{amsthm}
\usepackage{fancyhdr}
\usepackage{latexsym}
\usepackage{mathrsfs}
\usepackage{wasysym}
\usepackage{fancyhdr}
\usepackage{float}
\usepackage{graphicx}
\usepackage[numbers,sort&compress]{natbib}
\usepackage{setspace}
\usepackage{calligra}
\usepackage{xcolor} 

\allowdisplaybreaks[4]

\newtheorem{theorem}{\bf Theorem}[section]
\newtheorem{lemma}[theorem]{Lemma}
\newtheorem{proposition}[theorem]{Proposition}

\newtheorem{problem}[theorem]{Problem}
\newtheorem{corollary}[theorem]{Corollary}
\newtheorem{conj}[theorem]{Conjecture}

\begin{document}
\begin{spacing}{1.1}
\title{New results on openly disjoint cycles in regular digraphs}
\author{Bin Chen$^a$, \, Wenshen Wu$^a$, \, Meishuang Chen$^b$\thanks {Email addresses: cmsfzu07@163.com}\\
\small$^a$School of Mathematics and Statistics, Fuzhou University, Fujian, China\\
\small $^b$Center for Discrete Mathematics, Fuzhou University, Fujian, China}
\date{}
\maketitle
{\bf Abstract}:
Given a digraph \(D\), let \(c(D)\) be the largest size of a collection of
openly disjoint directed cycles through a common vertex, and let \(c_r\) be the
minimum of \(c(D)\) over all \(r\)-regular digraphs \(D\). In 2005, Seymour
asked whether every \(r\)-regular digraph satisfies \(c(D)=r\). Mader answered Seymour's
problem negatively for every \(r\ge 8\), but conjectured that nevertheless
\(\lim_{r\to\infty}c_r=\infty\). Recently, Steiner confirmed this conjecture by
proving
\[
\left\lceil \frac{3r}{22}\right\rceil
\le c_r \le
7\left\lceil \frac{r}{8}\right\rceil.
\]

In this paper, we first answer Seymour's
problem negatively for $r\in \{6,7\}$, and then improve Steiner's bound to 
\[
\left\lceil \frac{r}{6}\right\rceil
\le c_r \le
5\left\lceil \frac{r}{6}\right\rceil.
\]

{\bf AMS}\,: 05C20; 05C38.


\section{Introduction}

\noindent

All digraphs considered in this paper are simple, that is, they contain neither loops nor parallel arcs. Directed cycles\;(cycles for short) are one of the most fundamental objects in digraph theory. Given a digraph $D$, a collection of cycles \(C_1,C_2,\ldots,C_k\) containing a common vertex \(v\) is called \emph{openly disjoint} if
\[
V(C_i)\cap V(C_j)=\{v\}
\qquad\text{for all } i,j\in \{1,2,\ldots,k\}\; \text{with}\;i\neq j.
\]
For any vertex $v\in V(D)$, let \(c_v(D)\) be the maximum size of a collection of openly disjoint cycles
through \(v\), and we denote
\[
c(D)=\max_{v\in V(D)} c_v(D).
\]

Scholars are interested in finding openly disjoint cycles under the condition of minimum out-degree. In \cite{Thomassen1985a}, Thomassen proved that every digraph $D$ with minimum out-degree at least two satisfies that $c(D)\geq 2$. Naturally, we wonder whether the condition that the minimum out-degree is large enough can guarantee the existence of a large number of openly disjoint cycles through a common vertex. Surprisingly, Thomassen \cite{Thomassen1985b} answered this problem in the negative by constructing a family of digraphs in which every digraph $D$ has arbitrarily large minimum out-degree and minimum in-degree but $c(D)\leq 2$. Indeed, the digraphs constructed by Thomassen are far from being regular.

In 2005, Seymour proposed the following problem regarding the existence of many openly disjoint cycles in regular digraphs, where an \emph{$r$-regular digraph} is a digraph in which every vertex has both out-degree and in-degree $r$. 

\begin{problem}[Seymour \cite{Mader2010}] \label{seymour}
Does every \(r\)-regular digraph \(D\) satisfy that $c(D)=r$?
\end{problem}

It is worth mentioning that Seymour's problem is highly relevant to the conjecture of Behzad, Chartrand and Wall
from 1970 \cite{BCW1970}, which demonstrates that every \(r\)-regular digraph with
girth at least \(g\) has at least \(r(g-1)+1\) vertices. Actually, if
Problem~\ref{seymour} is true, then the \(r\) openly disjoint
cycles through a vertex \(x\) would each have length at least \(g\), and hence they would together contain at least
\(1+r(g-1)\) vertices, implying the conjecture of Behzad, Chartrand and Wall.
There are also two stronger conjectures, one is the Caccetta-Häggkvist
conjecture \cite{CH1978}, which states that every digraph with minimum out-degree at least $r$ and girth at least $g$ contains at least $r(g-1)+1$ vertices,  and the other is Hoàng-Reed conjecture \cite{HR1987}, which asserts that every digraph with minimum out-degree at least $r$ contains a collection of cycles $C_1,C_2,\ldots,C_r$ such that
$
\left|
V(C_i)\cap\bigcup_{j=1}^{i-1}V(C_j)
\right|\leq 1$ for every $2\leq i\leq r.$ Despite extensive efforts, these conjectures remain open in general. We refer the reader to \cite{Bondy,Shen,HSA,HHK,HKN,Wel} for partial results, and to \cite{Sullivan} for a survey of further related results and problems.

We will introduce the parameter $c_r$ which has appeared in \cite{Mader2010, Steiner2026}. Let $c_r$ be the minimum value of $c(D)$ among all $r$-regular digraphs. Problem~\ref{seymour} asks whether $c_r=r$. Observe that, for any positive integer $r$, $c_r\leq r$ , and $c_1=1$. As mentioned earlier, Thomassen \cite{Thomassen1985a} verified that $c_2=2$. Later in 2009, Mader \cite{Mader2010} obtained an affirmative result for the case $r=3$, that is, $c_3=3$. In the same paper, Mader also gave positive
answers for vertex-transitive digraphs and for digraphs with strong
vertex-connectivity at least \(r-2\). However, Mader proved that $c_r\leq r-1$ for every \(r\geq 8\) by  
constructing an \(r\)-regular digraph \(D\) with \(c(D)\leq r-1\). This indicates that Problem~\ref{seymour} has a negative answer for every $r\geq 8$. 
Prior to the present work, the cases
$r\in\{4,5,6,7\}$ remained open. In \cite{Mader2010-2}, Mader  showed that \(c_r\geq 4\) for every \(r\geq 7\).

Mader \cite{Mader2010} also constructed a $2k^2$-regular digraph $D$ for any integer $k\geq 2$ with $c(D)\leq 2k^2-k+1$, which implies that \(c_r\leq r-\Theta(\sqrt{r})\). Furthermore, Mader in \cite{Mader2010} conjectured that \(c_r\) should grow arbitrarily large with
\(r\).

\begin{conj}[{\cite{Mader2010}}] \label{mader}
\(\displaystyle \lim_{r\to\infty} c_r = \infty\).
\end{conj}

Recently, Steiner \cite{Steiner2026} confirmed
this conjecture in a strong form.

\begin{theorem}[\cite{Steiner2026}] \label{steiner}
For any $r\in \mathbb{N}$, $$\lceil \frac{3r}{22}\rceil \leq c_r\leq 7\lceil \frac{r}{8}\rceil.$$
\end{theorem}

Steiner also proved that $c_r/r$ converges, yielding the following corollary.
\begin{corollary}[\cite{Steiner2026}] \label{coro}
The limit $\displaystyle \lim_{r\to\infty} \frac{c_r}{r}$ exists and lies between $\frac{3}{22}$ and $\frac{7}{8}$.
\end{corollary}

The same author further raised the following problem.
\begin{problem}[\cite{Steiner2026}] \label{Ste}
Determine  $\displaystyle \lim_{r\to\infty} \frac{c_r}{r}$.
\end{problem}

Motivated by these results, we shall further explore Problems~\ref{seymour} and~\ref{Ste}. Firstly, we will construct a $6$-regular\;(resp., $7$-regular) digraph $D$ satisfying that $c(D)\leq 5$\;(resp., $c(D)\leq 6$), that is, we provide counterexamples to Seymour's problem for $r\in \{6,7\}$ and show that $c_r\leq r-1$ for $r\in \{6,7\}$. 
By combining a result of Steiner \cite{Steiner2026}, we can deduce that $c_r\leq 5\lceil \frac{r}{6}\rceil$, reducing the leading coefficient of the upper bound in Theorem~\ref{steiner} from $7/8$ to $5/6$. Then we further improve the lower bound $\lceil \frac{3r}{22}\rceil$ in Theorem~\ref{steiner}  to $\lceil \frac{r}{6}\rceil$. 

\begin{theorem} \label{main}
For any $r\in \mathbb{N}$, $$\lceil \frac{r}{6}\rceil \leq c_r\leq 5\lceil \frac{r}{6}\rceil.$$
\end{theorem}

We thus obtain a straightforward corollary as below.
\begin{corollary}\label{maincoro}
The limit $\displaystyle \lim_{r\to\infty}\frac{c_r}{r}$ exists and lies between $\frac{1}{6}$ and $\frac{5}{6}$.
\end{corollary}

The paper is organized as follows. In Section 2, we give some notation and several auxiliary
results. In Section 3, we will present our main results, including the constructions of counterexamples to Problem~\ref{seymour} for $r\in\{6,7\}$, the improved lower bound on $c_r$, and the proofs of Theorem~\ref{main} and Corollary~\ref{maincoro}. We provide some remarks and discussions in the last
section.

\section{Preliminaries}

\noindent

For two positive integers $m\le n$, let $[m,n]=\{m, m+1, \ldots, n\}$, and $[n]=[1,n]$. 
Let $D=(V,A)$ be a digraph with vertex set $V=V(D)$ and arc set $A=A(D)$. For any arc $(u,v)$ of $D$, we call $u$ its \emph{tail} and $v$ its \emph{head}, and we say $u$ dominates $v$, or equivalently, $v$ is dominated by $u$. For simplicity, we write $u\rightarrow v$ for this arc $(u,v)$. In addition, we denote by $a(D)=|A(D)|$.

Let $v$ be a vertex of $D$. The set $N^{+}_{D}(v)=\{u\in V(D):v\rightarrow u\}$  is the  \emph{out-neighborhood} of $v$, and let $d^{+}_{D}(v)=|N^{+}_{D}(v)|$ be the \emph{out-degree} of $v$. The set $N^{-}_{D}(v)=\{u\in V(D):u\rightarrow v\}$ is the \emph{in-neighborhood} of $v$, and let $d^{-}_{D}(v)=|N^{-}_{D}(v)|$ be the \emph{in-degree} of $v$. The vertices in $N^{+}_{D}(v)$ and $N^{-}_{D}(v)$ are called the \emph{out-neighbors} and \emph{in-neighbors} of $v$, respectively. The \emph{maximum out-degree}\;(resp., \emph{minimum out-degree}) of $D$ is $\Delta^{+}(D)=\text{max}\{ d^{+}_{D}(v):v\in V(D)\}$\;(resp., $\delta^{+}(D)=\text{min}\{d^{+}_{D}(v):v\in V(D)\}$). Analogously, we use $\Delta^{-}(D)$ to mean the \emph{maximum in-degree} and $\delta^{-}(D)$ to mean the \emph{minimum in-degree} of $D$, respectively. A digraph $D$ is \emph{$r$-out-regular} if every vertex has out-degree $r$, and similarly a digraph $D$ is \emph{$r$-in-regular} if every vertex has in-degree $r$. A digraph is \emph{$r$-regular} if it is both $r$-out-regular and $r$-in-regular. An \emph{$(\alpha, \beta)$-vertex} $w$ of $D$ is a vertex with $d^+_D(w)=\alpha$ and $d^-_D(w)=\beta$. Observe that $c_w(D)\leq \text{min}\{\alpha,\beta\}$.

A digraph $H$ is a \emph{subdigraph} of $D$ if $V(H)\subseteq V(D)$ and $A(H)\subseteq A(D)$. For every subdigraph $H$ of $D$ and every
$v\in V(H)$, let
$
N_H^+(v)=\{u\in V(H):v\rightarrow u
\text{ in }H\}
$
and $d_H^+(v)=|N_H^+(v)|$.  For any $X\subseteq V(D)$, let $D[X]$ be the subdigraph induced by $X$ in $D$, and let $D-X$ be the subdigraph $D[V(D)\backslash X]$. For every pair of disjoint vertex sets $X$ and $Y$, we denote $A(X,Y)$ to be those arcs between $X$ and $Y$ with tails in $X$ and heads in $Y$, and write $a_D(X,Y)$ for the cardinality of $A(X,Y)$. Moreover, we denote by $a_D[X,Y]=a_D[Y,X]=a_D(X,Y)+a_D(Y,X)$ the number of arcs $(u,v)$ between $X$ and $Y$ with $u\in X$ and $v\in Y$ or $u\in Y$ and $v\in X$.  Additionally, a set $X\subseteq V(D)$ is said to be an \emph{independent set} if there is no arc in $D[X]$.

Next, we shall provide some auxiliary results. The first result tells us that the size of every digraph $D$ can be bounded by a function of its order and $c(D)$.

\begin{lemma}[\cite{Steiner2026}]  \label{lemma 1}
Every digraph $D$ on $n$ vertices satisfies that $a(D)<\frac{n(n+c(D))}{2}$.
\end{lemma}

For any digraph $D$, a partition  $(X,Y)$ of $V(D)$ consists of $X\subseteq V(D)$ and $Y\subseteq V(D)$ with $X\cup Y=V(D)$ and $X\cap Y=\emptyset$. We also need the following structural result.

\begin{lemma}[\cite{Steiner2026}]  \label{lemma 2}
Let $D$ be a digraph. Then for every vertex $v\in V(D)$, there exists a set $S\subseteq V(D)$ with $|S|\leq c(D)$ and a partition $(X,Y)$ of $V(D)\backslash S$ such that $N^+_D(v)\subseteq X\cup S$, $N^-_D(v)\subseteq Y\cup S$ and $A(X,Y)=\emptyset$.
\end{lemma}

Then we present the following lemma, which can estimate the upper bound.
\begin{lemma}[\cite{Steiner2026}] \label{lemma 3}
For every $r,m\in \mathbb{N}$, $c_{rm}\leq c_{r}*m$.
\end{lemma}

We end this section with another auxiliary result as below, which will play an important role in the improvement of the lower bound.
\begin{lemma}\label{lemma 4}
Let $r$ be a positive integer and let $D$ be a digraph on $n\geq r/2$ vertices. If
$a(D)\geq rn-\frac{r^2}{3}$,
$\Delta^+(D),\Delta^-(D)\leq r$
and
$
c(D)<r/6,
$
then
$n>4r/3,$
and there exists a vertex $v\in V(D)$ such that both $d_D^+(v)$ and $d_D^-(v)$ are at least $r/2$.
\end{lemma}
\noindent\textbf{Proof of Lemma~\ref{lemma 4}.} Set $x=n/r$, and it is clear that $x\geq 1/2$ since $n\geq r/2$.
By Lemma~\ref{lemma 1} and the conditions $a(D)\geq rn-\frac{r^2}{3}$ and $c(D)<r/6$, we have that
\[
rn-\frac{r^2}{3}
\leq a(D)
<
\frac{n(n+c(D))}{2}
<
\frac{n^2}{2}+\frac{rn}{12}.
\]
Dividing by $r^2$ and rearranging yields that
\[
x^2-\frac{11}{6}x+\frac{2}{3}>0.
\]
Equivalently,
\[
\left(x-\frac{1}{2}\right)
\left(x-\frac{4}{3}\right)>0.
\]
Since $x\geq 1/2$, it follows that $x>4/3$, and hence $n>4r/3$.

Next, we will verify the second statement. Suppose to the contrary that no vertex has both out-degree and in-degree
at least $r/2$.
Then, for every vertex $v\in V(D)$, at least one of
$d_D^+(v)$ and $d_D^-(v)$ is less than $r/2$.
As $d_D^+(v)\leq \Delta^+(D)\leq r$ and  $d_D^-(v)\leq \Delta^-(D)\leq r$, we thus get
$d_D^+(v)+d_D^-(v)<\frac{3r}{2}$.
Summing over all vertices, one can deduce that
\[
2a(D)
=
\sum_{v\in V(D)}
\bigl(d_D^+(v)+d_D^-(v)\bigr)
<
\frac{3rn}{2}.
\]
Combining this with the assumed lower bound on $a(D)$, we conclude that
\[
rn-\frac{r^2}{3}<\frac{3rn}{4},
\]
which yields that $n<4r/3$, contradicting the first statement that $n>4r/3$. This proves that there indeed exists a vertex of both 
out-degree and in-degree at least $r/2$. \hfill $\blacksquare$

\section{Main results}

\noindent

In this section, we establish our main results. In Subsection~\ref{r=6,7}, we construct a $6$-regular digraph $D$ with $c(D)\leq 5$ and a $7$-regular digraph $D$ with $c(D)\leq 6$, thereby providing counterexamples to Problem~\ref{seymour}. In Subsection~\ref{improvementlower}, we improve the lower bound on $c_r$ to $\left\lceil\frac{r}{6}\right\rceil$. The proofs of Theorem~\ref{main} and  Corollary~\ref{maincoro} will be presented in Subsection~\ref{thm-coro}.

\subsection{Counterexamples to Problem~\ref{seymour} for $r\in\{6,7\}$}\label{r=6,7}

\noindent

\begin{figure}[!hbpt]
\begin{center}
\includegraphics[scale=0.31]{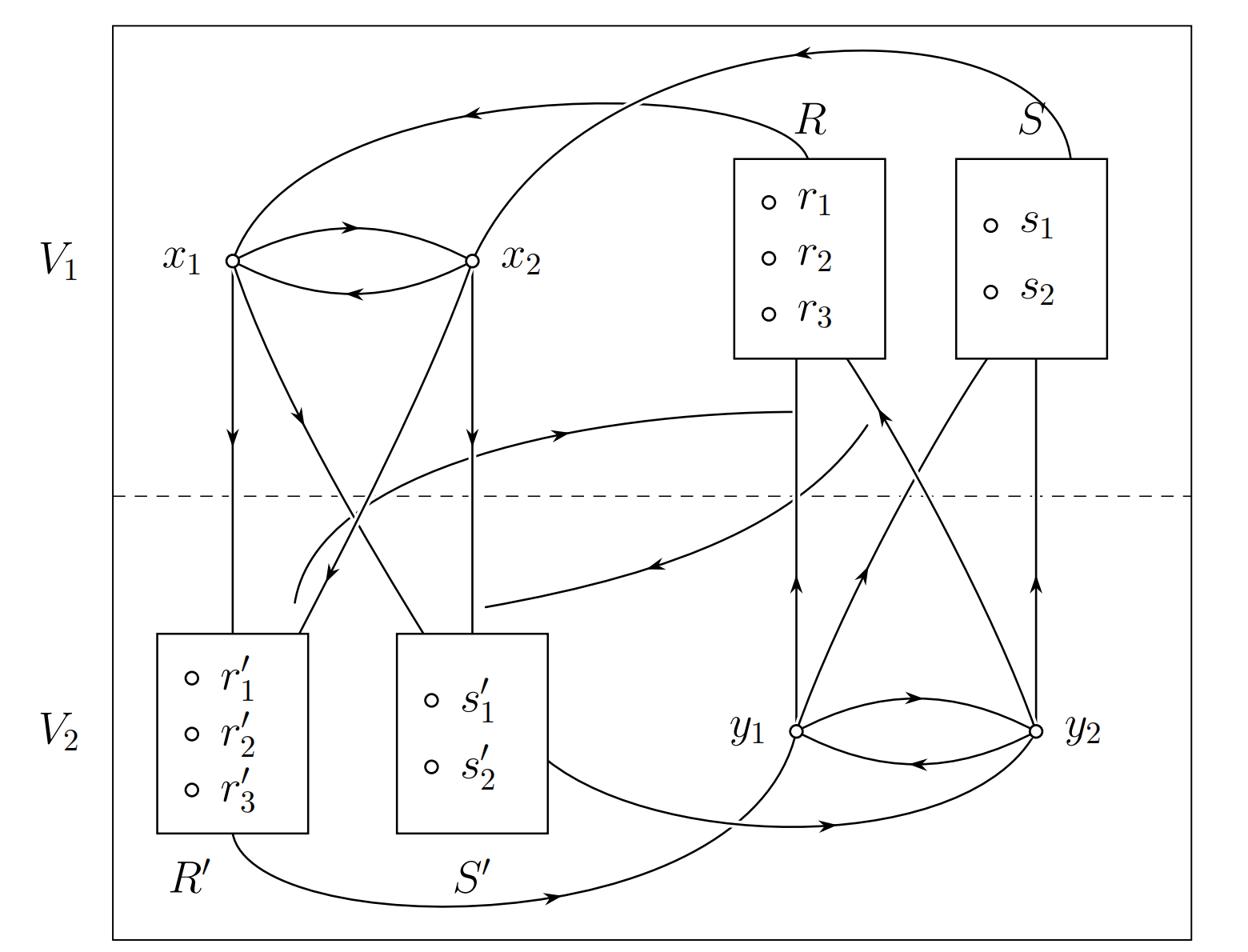}
  \caption{\small
   $F$.
  }
  \label{F}
\end{center}
\end{figure}

In this subsection, we would like to give counterexamples to Problem~\ref{seymour} for $r\in\{6,7\}$. Firstly, we will construct a $7$-regular digraph $H_1$ satisfying $c(H_1)\leq 6$. The construction of $H_1$ relies on the following digraph $F$\;(see Figure~\ref{F}).

\vspace{0.2cm}

\noindent\textbf{Construction of $F$.} Let $F$ be a digraph with vertex set $V(F)=\{x_1,x_2,y_1,y_2\}\cup (R\cup S)\cup (R'\cup S')$, where $R=\{r_1,r_2,r_3\}$ and $S=\{s_1,s_2\}$, and $R'=\{r_1',r_2',r_3'\}$ and $S'=\{s_1',s_2'\}$, and arc set
$A(F)$ consisting of $x_1\to x_2$, $x_2\to x_1$, $y_1\to y_2$, $y_2\to y_1$, $r_i\to x_1$ for each $i\in[3]$, $s_i\to x_2$ for each $i\in[2]$, $r_i'\to y_1$ for each $i\in[3]$, $s_i'\to y_2$ for each $i\in[2]$, together with all arcs from
$\{x_1,x_2\}\cup R\cup S$ to $R'\cup S'$
and all arcs from
$\{y_1,y_2\}\cup R'\cup S'$ to $R\cup S$.

Notice that $F$ is a $6$-out-regular digraph on $14$ vertices. Moreover, $x_1$ and $y_1$ are (6,4)-vertices, $x_2$ and $y_2$ are (6,3)-vertices, and every vertex in $R\cup S\cup R'\cup S'$ is a $(6,7)$-vertex. Next, we will show that $c(F)=5$.

 \begin{proposition}\label{exact5}
We have $c(F)=5$.
\end{proposition}
\noindent\textbf{Proof of Proposition~\ref{exact5}.}  For each $w\in R'\cup S'$, the arcs $r_1\to w$ and $w\to r_1$ form a cycle of length two. These five cycles are openly disjoint, so $c(F)\geq c_{r_1}(F)\geq 5$.
 As every $(\alpha,\beta)$-vertex $v$ satisfies that $c_v(D)\leq \text{min}\{\alpha,\beta\}$, we thus have that, for any vertex $z\in \{x_1,x_2,y_1,y_2\}$, $c_z(F)<5$. It remains to consider vertices $z\in (R\cup S)\cup (R'\cup S')$. By symmetry, we only need to take the case that $z\in R\cup S$. When $z=r_i$ for some $i\in[3]$, we will show that $c_z(F)\leq 5$. Assume towards a contradiction that $c_z(F)\geq 6$. This means that there are $6$ openly disjoint cycles passing through $z$. The 6 cycles leave $z$ through distinct out-neighbors. Since $N_F^+(z)=\{x_1\}\cup R'\cup S'$, one of them uses the arc $z\to x_1$, while the other five leave $z$ through the five distinct
vertices of $R'\cup S'$. Note that every directed path with initial vertex $x_1$ and terminal vertex $z$ must contain a vertex in $R'\cup S'$. In other words, any cycle containing the arc $(z,x_1)$ must contain a vertex in $R'\cup S'$. Consequently, at least two of these cycles share a vertex
in $R'\cup S'$, contradicting the assumption that
they are openly disjoint. Hence, $c_z(F)\leq 5$. By using an almost identical method, one can also prove that $c_{z}(F)\leq 5$ when $z=s_i$ for $i\in \{1,2\}$. Therefore, we conclude that $c(F)=5$.    
\hfill $\blacksquare$

\vspace{0.2cm}

For convenience, we denote $X=\{x_1,x_2\}$, $Y=\{y_1,y_2\}$, $V_1=X\cup R\cup S$ and $V_2=Y\cup R'\cup S'$. Note that $V(F)=V_1\cup V_2$ and $V_1\cap V_2=\emptyset$. Let $F_1$, $F_2$, $F_3$ and $F_4$ be pairwise disjoint copies of $F$, and let $\overleftarrow{F}_i$ be the digraph obtained from $F_i$ by reversing all of its arcs, where $i\in[4]$.  In addition, let $B$ and $B'$ be disjoint independent sets of size two, with vertices
$b_1,b_2$ and $b'_1,b'_2$, respectively, disjoint from $V(F_1)\cup V(F_2)\cup V(F_3)\cup V(F_4)$.
For each $i\in[4]$, let $R_i,R_i',S_i,S_i',X_i,Y_i$
denote the copies of $R,R',S,S',X,Y$ in $F_i$, respectively. For any $i\in[4]$ and $j\in[2]$, let
$V_{i,j}$, $x_{i,j}$ and $y_{i,j}$ denote the copies of $V_j$, $x_j$ and $y_j$ in $F_i$, respectively.

\begin{figure}[!hbpt]
\begin{center}
\includegraphics[scale=0.31]{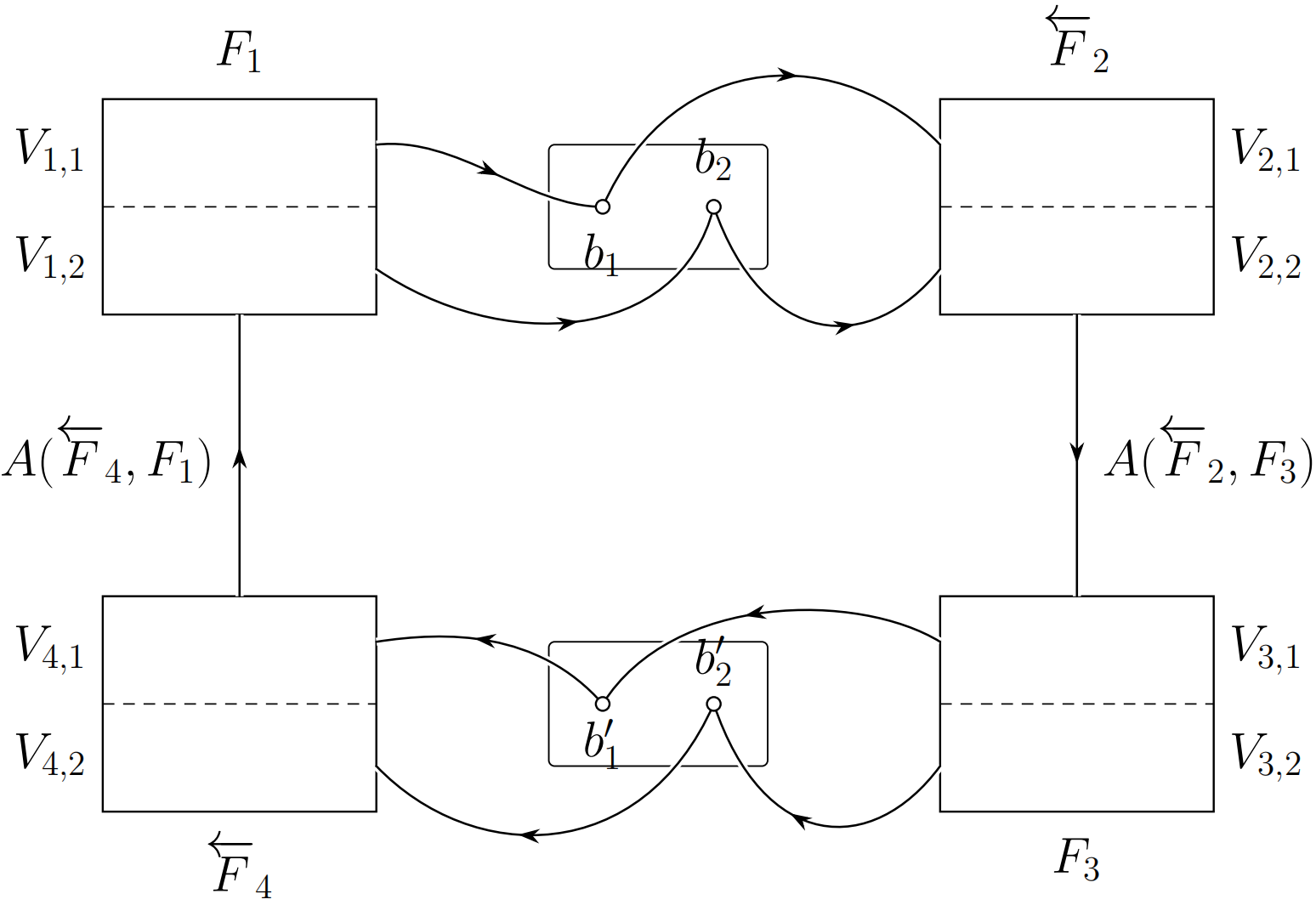} 
  \caption{\small
    $H_1$.
  }
  \label{H_1}
\end{center}
\end{figure}

Now we are ready to construct a $7$-regular digraph $H_1$ as below\;(see Figure~\ref{H_1}).

\noindent\textbf{Construction of $H_1$.} Let $H_1$ be a digraph obtained from the disjoint union
$F_1\cup\overleftarrow{F}_2\cup F_3\cup\overleftarrow{F}_4\cup B\cup B'$
by adding all arcs from 
$V_{1,1}$ to $b_1$, from $V_{1,2}$ to $b_2$,
from $b_1$ to $V_{2,1}$, from $b_2$ to $V_{2,2}$,
from $V_{3,1}$ to $b'_1$, from $V_{3,2}$ to $b'_2$,
from $b'_1$ to $V_{4,1}$, from  $b'_2$ to $V_{4,2}$, and all arcs from $X_4\cup Y_4$ to $X_1\cup Y_1$ except two arcs $(y_{4,1}, x_{1,1})$ and $(x_{4,1}, y_{1,1})$,
from $X_2\cup Y_2$ to $X_3\cup Y_3$ except two arcs $(x_{2,1}, y_{3,1})$ and $(y_{2,1}, x_{3,1})$.

\vspace{0.2cm}

It is not difficult to check that $H_1$ is a $7$-regular digraph on $60$ vertices. Let us now show that $c(H_1)=6$. 

 \begin{proposition}\label{exact6}
We have $c(H_1)=6$.
\end{proposition}
\noindent\textbf{Proof of Proposition~\ref{exact6}.} Let us first show that $c(H_1)\geq 6$. Let $z=s_{1,1}\in S_1$ be the vertex corresponding
to $s_1$ in $F_1$. We will state that $c_{s_{1,1}}(H_1)\geq 6$. Indeed, by the construction of $H_1$, there must be a path starting at $s_{1,1}$ and terminating at $y_{4,1}$ which is disjoint from $R_1'\cup S_1'$. Together with  $(y_{4,1},y_{1,1})$ and  $(y_{1,1},s_{1,1})$, we can thus find a cycle through $s_{1,1}$, which is openly disjoint from 5 cycles contained in $H_1[(R_1'\cup S_1')\cup \{s_{1,1}\}]$. It follows that $c_{s_{1,1}}(H_1)\geq 6$, and hence $c(H_1)\geq 6$.

Next, we would like to verify that $c(H_1)\leq 6$. Firstly, we will prove that $c_{z}(H_1)\leq 6$ for any $z\in V(F_1)\cup V(\overleftarrow{F}_4)\cup B$. 
From the construction of $H_1$ we know that every cycle through a vertex of $B$ must contain some vertex in $B'$. This indicates that $c_{z}(H_1)\leq 2$ when $z\in B$. It remains to consider that $z\in V(F_1)\cup V(\overleftarrow{F}_4)$. 
We first consider the case that $z\in V(F_1)$. Note that $c_z(H_1)\leq c_z(F_1)+2$ by the construction of $H_1$ because every cycle through $z$ but not contained in $F_1$ must pass through some vertex of $B$ and $|B|=2$. Additionally, we know that $x_{1,1}$ and $y_{1,1}$ are $(6,4)$-vertices of $F_1$, and    $x_{1,2}$ and $y_{1,2}$ are $(6,3)$-vertices of $F_1$. Thereby, if $z\in X_1\cup Y_1$, then $c_z(F_1)\leq 4$, yielding that $c_z(H_1)\leq 6$. We next assume that $z\in (R_1\cup S_1)\cup (R_1'\cup S_1')$. Analogously, we only need to take the case that $z\in R_1\cup S_1$ as the other case is very similar. By Proposition~\ref{exact5}, we have that $c(F)=5$. From the previous results, we can easily deduce that $c_z(H_1)\leq 7$. We would like to show that $c_z(H_1)\leq 6$. Suppose, on the contrary, that $c_z(H_1)=7$. Consider 7 openly disjoint cycles through $z$. At most 5 are contained in $F_1$, and at most two are not. Hence exactly five are contained in $F_1$. The remaining two cycles each contain exactly one vertex of \(B\), with one containing \(b_1\) and the other containing \(b_2\). Every cycle through $z$ contained in $F_1$ meets $R_1'\cup S_1'$. Since the five internal cycles are openly disjoint and $|R_1'\cup S_1'|=5$, together they contain every vertex of $R_1'\cup S_1'$.  Moreover, every directed path from $z$ to $b_2$ must pass through $R_1'\cup S_1'$ or $b_1$. Hence the cycle containing $b_2$ shares a vertex other than $z$ with one of the other six cycles, which leads to a contradiction. This proves that any vertex $z\in R_1\cup S_1$ satisfies that $c_z(H_1)\leq 6$. By using almost identical arguments, one can also infer that  $c_z(H_1)\leq 6$ when $z\in R_1'\cup S_1'$. As a consequence, we conclude that any $z\in V(F_1)$ satisfies that $c_z(H_1)\leq 6$.

Then we turn to the case $z\in V(\overleftarrow{F}_4)$. Every cycle through $z$ that is not contained in
$\overleftarrow{F}_4$ meets $B'$. Since $|B'|=2$, any collection of openly disjoint cycles through $z$ contains at most two such cycles. Hence $c_z(H_1)\leq c_z(\overleftarrow{F}_4)+2.$ As $x_{4,1}$ and $y_{4,1}$ are $(4,6)$-vertices of $\overleftarrow{F}_4$, and $x_{4,2}$ and $y_{4,2}$ are $(3,6)$-vertices, we have $ c_z(\overleftarrow{F}_4) \leq d_{\overleftarrow{F}_4}^+(z)\leq 4$ for every $z\in X_4\cup Y_4.$
Therefore, $c_z(H_1)\leq 6$ for every $z\in X_4\cup Y_4$.  Let us next assume that $z\in (R_4\cup S_4)\cup (R_4'\cup S_4')$. Similarly, it suffices to consider $z\in R_4\cup S_4$. Observe that $c(\overleftarrow{F}_4)=5$ and $c_z(H_1)\leq 7$. Our goal is to verify that $c_z(H_1)\leq 6$. Suppose towards a contradiction that $c_z(H_1)=7$. Consider 7 openly disjoint cycles through $z$.
At most five are contained in $\overleftarrow{F}_4$,
and at most two are not.
Hence exactly five are contained in $\overleftarrow{F}_4$.
The remaining two cycles contain exactly one vertex
of $B'$ each, one containing $b_1'$ and the other
containing $b_2'$. Every cycle through $z$ contained in $\overleftarrow{F}_4$ meets $R_4'\cup S_4'$. Since the five internal cycles are openly disjoint,
together they contain every vertex of $R_4'\cup S_4'$. Moreover, every arc with tail outside $V_{4,1}$
and head in $V_{4,1}$ has its tail in $(R_4'\cup S_4')\cup\{b_1'\}$. Thus every directed path from $b_2'$ to $z$
must pass through $R_4'\cup S_4'$ or $b_1'$. Hence the cycle containing $b_2'$ shares a vertex other than $z$ with one of the other six cycles, contradicting their being openly disjoint.  This proves that any vertex $z\in R_4\cup S_4$ satisfies that $c_z(H_1)\leq 6$. Also, we can deduce, by an analogous method, that  $c_z(H_1)\leq 6$ when $z\in R_4'\cup S_4'$. Consequently, we derive that any $z\in V(\overleftarrow{F}_4)$ satisfies that  $c_z(H_1)\leq 6$.

Symmetrically, by applying nearly identical arguments, we can also deduce that $c_{z}(H_1)\leq 6$ for $z\in V(\overleftarrow{F}_2)\cup V(F_3) \cup B'$. In summary, we obtain that any vertex $z\in V(H_1)$ satisfies that $c_z(H_1)\leq 6$, and therefore $c(H_1)=6$.   \hfill $\blacksquare$

\begin{figure}[!hbpt]
\begin{center}
\includegraphics[scale=0.31]{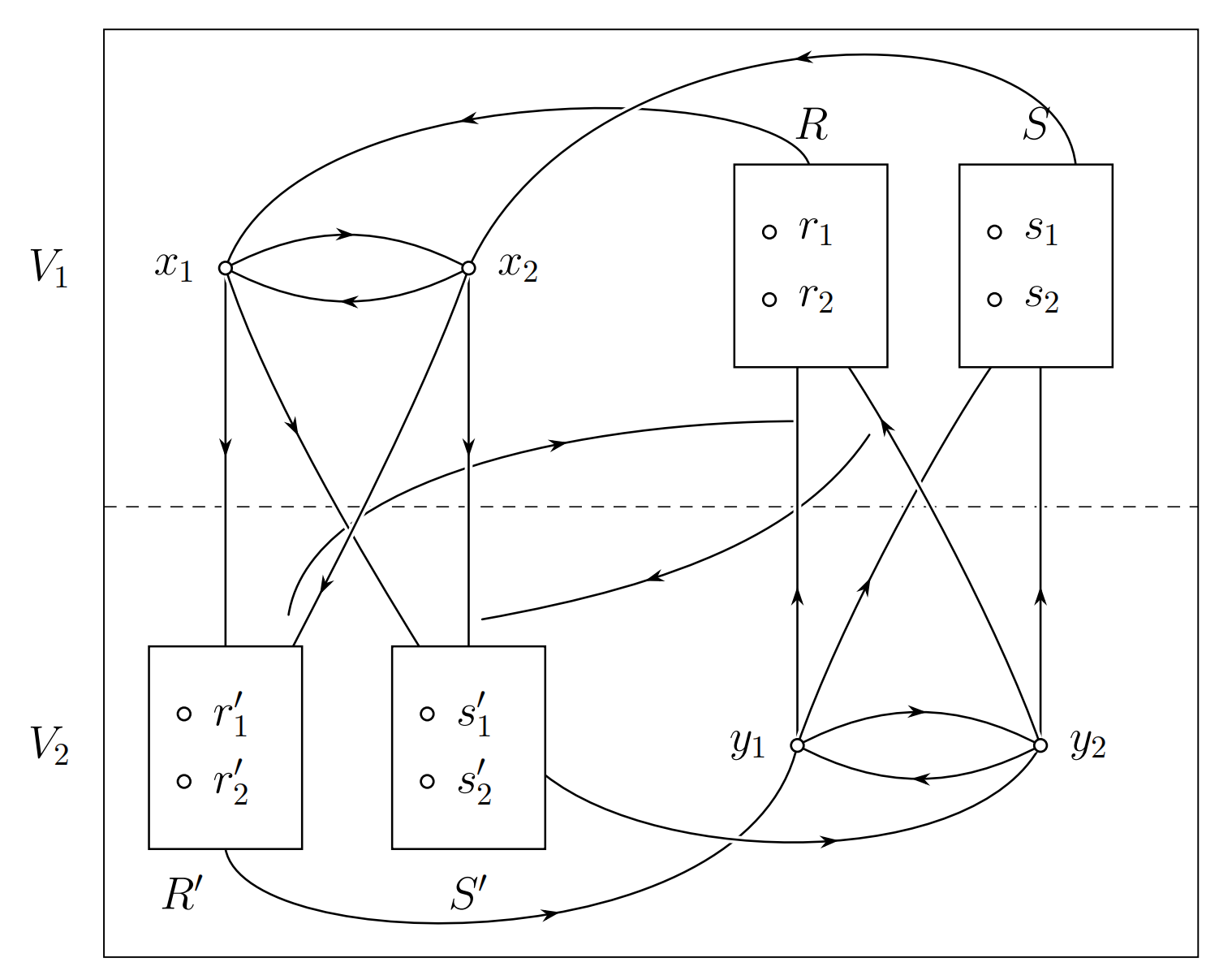}
  \caption{\small
    $G$.
  }
  \label{G}
\end{center}
\end{figure}

In the following, we will provide a $6$-regular digraph such that every vertex lies on at most 5 openly disjoint cycles. We demonstrate that the construction of this digraph will depend on a 5-out-regular digraph $G$\;(see Figure~\ref{G}).

\vspace{0.2cm}

\noindent\textbf{Construction of $G$.} Let $G$ be the digraph obtained from $F$
by deleting the vertices $r_3$ and $r_3'$.

\vspace{0.2cm}

Observe that $G$ is a $5$-out-regular digraph on $12$ vertices. It is evident that $x_1$, $x_2$, $y_1$ and $y_2$ are (5,3)-vertices, and every other vertex of $G$ is a $(5,6)$-vertex. We can obtain that $c(G)=4$ by the similar argument as Proposition~\ref{exact5}.

 \begin{proposition}\label{exact4}
We have $c(G)=4$.
\end{proposition}

\begin{figure}[!hbpt]
\begin{center}
\includegraphics[scale=0.31]{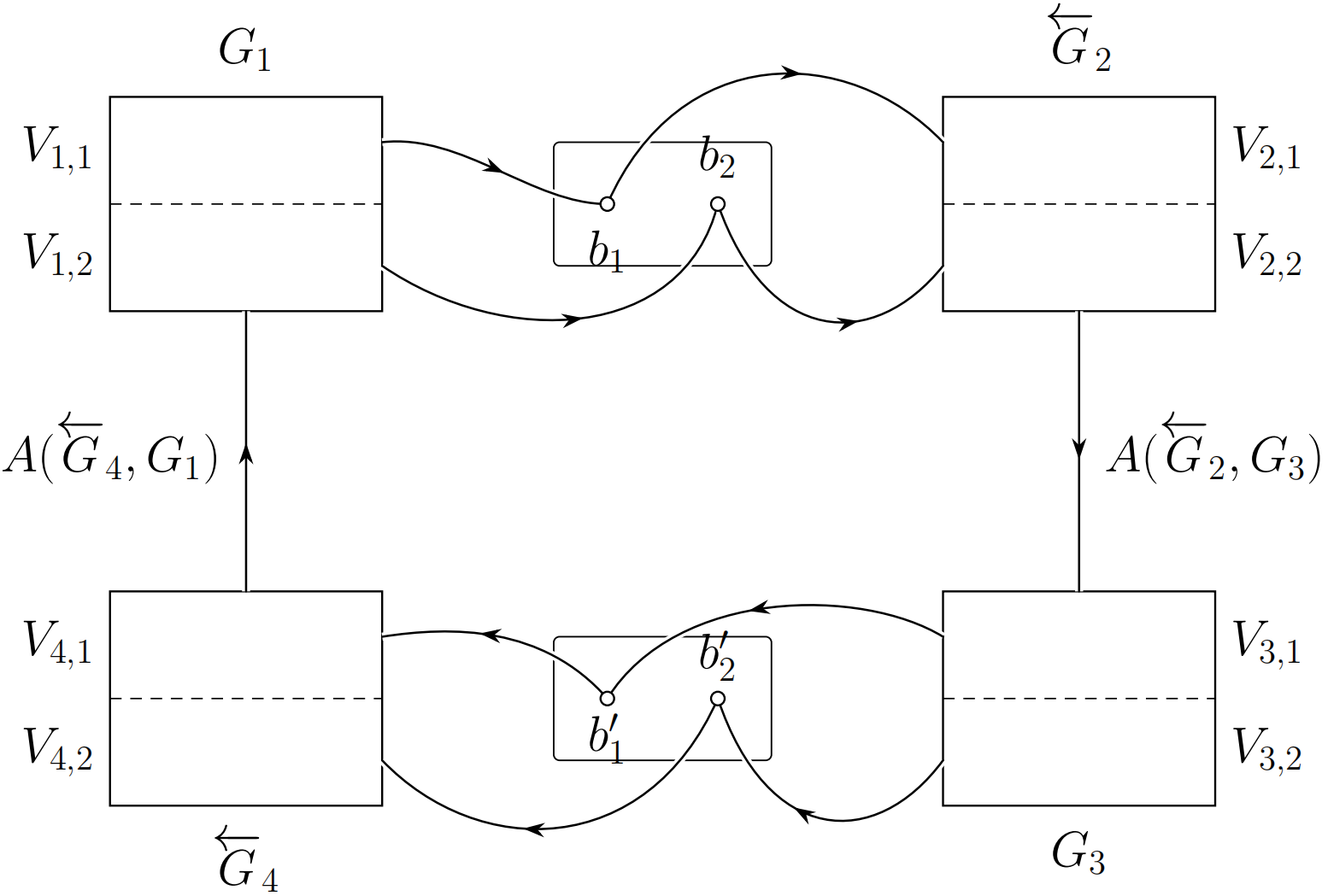}
  \caption{\small
    $H_2$.
  }
  \label{H_2}
\end{center}
\end{figure}

Moving on, we shall construct a $6$-regular digraph $H_2$ with $c(H_2)=5$\;(also see Figure~\ref{H_2}). We mark the vertices and vertex sets in $G$ to be consistent with those in $F$, and also let $B=\{b_1,b_2\}$\;(resp., $B'=\{b_1',b_2'\}$) be an independent set of size two. Let $G_1,G_2,G_3,G_4$ be pairwise vertex-disjoint copies of $G$, and let $\overleftarrow{G}_i$
be obtained from $G_i$ by reversing all its arcs. We use the corresponding notation for the vertex
subsets in these copies. In particular, $|V_{i,1}|=|V_{i,2}|=6$ for every $i\in[4]$. In this construction, $R_i$ and $R_i'$ denote the copies of $R\setminus\{r_3\}$ and
$R'\setminus\{r_3'\}$, respectively.

\vspace{0.2cm}

\noindent\textbf{Construction of $H_2$.} Let $H_2$ be a digraph obtained from the disjoint union
$G_1\cup\overleftarrow{G}_2\cup G_3\cup\overleftarrow{G}_4\cup B\cup B'$
by adding all arcs from 
$V_{1,1}$ to $b_1$, from $V_{1,2}$ to $b_2$,
from $b_1$ to $V_{2,1}$, from $b_2$ to $V_{2,2}$,
from $V_{3,1}$ to $b'_1$, from $V_{3,2}$ to $b'_2$,
from $b'_1$ to $V_{4,1}$, from  $b'_2$ to $V_{4,2}$, and all arcs from $X_4\cup Y_4$ to $X_1\cup Y_1$ except 4 arcs $(x_{4,1}, y_{1,2})$, $(x_{4,2}, y_{1,1})$, $(y_{4,1}, x_{1,2})$ and $(y_{4,2}, x_{1,1})$,
from $X_2\cup Y_2$ to $X_3\cup Y_3$ except 4 arcs $(x_{2,1}, y_{3,2})$, $(x_{2,2}, y_{3,1})$, $(y_{2,1}, x_{3,2})$ and $(y_{2,2}, x_{3,1})$.

\vspace{0.2cm}

One can check that $H_2$ is a $6$-regular digraph on $52$ vertices. By using almost identical method as that of Proposition~\ref{exact6}, one can derive that $c(H_2)=5$. 

 \begin{proposition}\label{exact}
We have $c(H_2)=5$.
\end{proposition}

From the above results, we conclude that there is a $7$-regular digraph $H_1$ with $c(H_1)\leq 6$ and there is also a  $6$-regular digraph $H_2$ with $c(H_2)\leq 5$(see Figure~\ref{H_2}). This means that $c_7\leq 6$ and $c_6\leq 5$.

\subsection{Improvement of the lower bound}\label{improvementlower}

\noindent

In this subsection, we aim to give a better lower bound on $c(D)$ for any $r$-regular digraph $D$.

\begin{theorem}\label{1/6}
For any $r\in \mathbb{N}$, every $r$-regular digraph $D$ on $n$ vertices
satisfies that
$c(D)\geq \left\lceil\frac{r}{6}\right\rceil$.
\end{theorem}

\noindent\textbf{Proof of Theorem~\ref{1/6}.} Suppose, on the contrary, that the assertion is not true. Then
$c(D)<\frac{r}{6}$.
Pick a subdigraph $D^*\subseteq D$ of minimum order, denoted by $n^*=|V(D^*)|$, subject to
\[
n^*\geq r/2
\qquad\text{and}\qquad
a(D^*)\geq rn^*-\frac{r^2}{3}.
\]
Observe that such a subdigraph must exist because $D$ itself satisfies
these conditions.

Since $D^*\subseteq D$, we have that
\[
\Delta^+(D^*),\Delta^-(D^*)\leq r
\qquad\text{and}\qquad
c(D^*)\leq c(D)<\frac{r}{6}.
\]
From Lemma~\ref{lemma 4} we know that
$n^*>\frac{4r}{3}$, and there is a vertex $v\in V(D^*)$ with
$d_{D^*}^+(v)\geq r/2$ and $d_{D^*}^-(v)\geq r/2$.

By applying Lemma~\ref{lemma 2} to $v$ in $D$, we can find a set $S\subseteq V(D)$ with
\[
|S|\leq c(D)<\frac{r}{6}
\]
and a partition $(X_S,Y_S)$ of $V(D)\setminus S$ such that
\[
N_D^+(v)\subseteq X_S\cup S,
\qquad
N_D^-(v)\subseteq Y_S\cup S,
\qquad
a_D(X_S,Y_S)=0.
\]
Additionally, we denote by 
\[
X_S^*=X_S\cap V(D^*),
\qquad
Y_S^*=Y_S\cap V(D^*),
\qquad
S^*=S\cap V(D^*).
\]
One easily sees that
$|X_S^*|+|Y_S^*|=n^*-|S^*|\geq n^*-|S|>
\frac{4r}{3}-\frac{r}{6}>r$.
This indicates that at least one of $X_S^*$ and $Y_S^*$ contains more than $r/2$ vertices.

Define a partition $(X^*,Y^*)$ of $V(D^*)$ by
\[
(X^*,Y^*)=
\begin{cases}
(X_S^*,\,Y_S^*\cup S^*), & \text{if}\; |X_S^*|>r/2,\\
(X_S^*\cup S^*,\,Y_S^*), & \text{otherwise}.
\end{cases}
\]
Notice that
\[
N_{D^*}^+(v)\subseteq X_S^*\cup S^*,
\qquad
N_{D^*}^-(v)\subseteq Y_S^*\cup S^*.
\]
Together with $d_{D^*}^+(v)\geq r/2$ and $d_{D^*}^-(v)\geq r/2$, one can infer that
$|X^*|,|Y^*|\geq r/2$ in either case.

Next, we will show that 
$$
a_{D^*}[X^*,Y^*]\geq \frac{r^2}{3}.
$$
Indeed, by the minimality of $D^*$, we can obtain that
\[
a(D^*[X^*])<r|X^*|-\frac{r^2}{3}
\qquad\text{and}\qquad
a(D^*[Y^*])<r|Y^*|-\frac{r^2}{3}.
\]
By combining $|X^*|+|Y^*|=|V(D^*)|=n^*$ and the fact $a(D^*)\geq rn^*-\frac{r^2}{3}$, we can thus get 
$$
a_{D^*}[X^*,Y^*]=a(D^*)-a(D^*[X^*])-a(D^*[Y^*])\geq rn^*-\frac{r^2}{3}-
\left(r|X^*|-\frac{r^2}{3}\right)
-\left(r|Y^*|-\frac{r^2}{3}\right)
\geq \frac{r^2}{3},
$$
as desired.

Since $D$ is $r$-regular, the numbers of arcs entering and leaving any vertex subset are equal. It follows that 
$a_D(X_S\cup S,Y_S)=a_D(Y_S,X_S\cup S)$ and $a_D(X_S,Y_S\cup S)=a_D(Y_S\cup S,X_S)$.

If $(X^*,Y^*)=(X_S^*,\,Y_S^*\cup S^*)$, then
$$a_{D^*}[X^*,Y^*]\leq a_D[X_S,Y_S\cup S]=2a_D(X_S,Y_S\cup S)=2a_D(X_S,S)
\leq 2r|S|<2r*\frac{r}{6}=\frac{r^2}{3},
$$
which contradicts that $a_{D^*}[X^*,Y^*]\geq \frac{r^2}{3}$.

If $(X^*,Y^*)=(X_S^*\cup S^*,Y_S^*)$, then 
$$a_{D^*}[X^*,Y^*]\leq a_D[X_S\cup S,Y_S]=2a_D(X_S\cup S,Y_S)=2a_D(S,Y_S)\leq 2r|S|
<2r*\frac{r}{6}
=\frac{r^2}{3},$$
which is also a contradiction.

In summary, we can obtain that $c(D)\geq r/6$.
As $c(D)$ is an integer, we conclude that
$c(D)\geq \left\lceil\frac{r}{6}\right\rceil$, completing our proof of Theorem~\ref{1/6}.  \hfill $\blacksquare$

\subsection{Proofs of Theorem~\ref{main} and Corollary~\ref{maincoro}}\label{thm-coro}

\noindent

In this subsection, we would like to provide the proofs of Theorem~\ref{main} and Corollary~\ref{maincoro}.

\vspace{0.2cm}

\noindent\textbf{Proof of Theorem~\ref{main}.} From Theorem~\ref{1/6} we know that $c_r\geq \left\lceil\frac{r}{6}\right\rceil$. It  remains to show that $c_r\leq 5\left\lceil\frac{r}{6}\right\rceil$. Indeed, as mentioned in Subsection~\ref{r=6,7}, there exists a $6$-regular digraph $D$ with $c(D)\leq 5$, which yields that $c_6\leq 5$. By Lemma~\ref{lemma 3}, one can see that $c_{6m}\leq 5m$ for any positive integer $m$. By combining the fact that $c_\alpha\leq c_\beta$ for any $\alpha,\beta$ with $\alpha\leq \beta$, we can thus deduce that $c_r\leq c_{6\left\lceil\frac{r}{6}\right\rceil}\leq  5\left\lceil\frac{r}{6}\right\rceil$, completing the proof of  Theorem~\ref{main}. \hfill $\blacksquare$

\vspace{0.2cm}

Let us now turn to the proof of Corollary~\ref{maincoro}. We clarify that this proof is almost identical to that of Steiner \cite{Steiner2026}, and we present the proof in detail to make this paper self-contained.

\vspace{0.2cm}

\noindent\textbf{Proof of Corollary~\ref{maincoro}.} We denote \(a_r=c_r/r\) and \(L=\inf_{r\in \mathbb{N}}a_r\). Our goal is to show that $a_r$ converges to $L$. To see this,  it is sufficient to show
\(\limsup_{r\to\infty}a_r\le L\) because 
\(\liminf_{r\to\infty}a_r\ge L\) is trivially true.

Fix a positive integer $r$. For every integer $n\geq r$, write $n=qr+s$, where $q\geq 1$ and $0\leq s<r$ are integers. By the monotonicity of $c_r$ and Lemma~\ref{lemma 3},
$c_n=c_{qr+s}\le c_{(q+1)r}\le (q+1)c_r$. Hence, we get
$
a_n=\frac{c_n}{n}
\le \frac{(q+1)c_r}{qr+s}
\le \frac{(q+1)c_r}{qr}$.
Letting \(n\to\infty\) (so \(q\to\infty\)) implies that
\[
\limsup_{n\to\infty}a_n\le \frac{c_r}{r}=a_r.
\]
Due to the arbitrariness of $r$, we thus have
\[
\limsup_{n\to\infty}a_n\le \inf_{r\ge1}a_r=L,
\]
which indicates that the limit \(L=\lim_{r\to\infty}c_r/r\) exists.

For the lower bound, \(c_r\ge \lceil \frac{r}{6}\rceil\ge r/6\) tells us that
\(a_r\ge 1/6\) for all \(r\), yielding that \(L\ge 1/6\). For the upper bound, due to
$c_r\le 5\lceil\frac{r}{6}\rceil\leq 5(\frac{r}{6}+1)=\frac{5r}{6}+5$, we can thus get $L\leq 5/6$.
Consequently, we conclude that $\frac{1}{6}\le \lim_{r\to\infty}\frac{c_r}{r}\le \frac{5}{6}$. \hfill $\blacksquare$

\section{Concluding remarks}

\noindent

In this paper, we mainly investigate Problems~\ref{seymour} and~\ref{Ste}, and obtain some new results. Problem~\ref{seymour} asks whether every $r$-regular digraph $D$ satisfies that $c(D)=r$, it has positive answers for $r=1$\;(trivially true), $r=2$\;(due to Thomassen \cite{Thomassen1985a}) and $r=3$\;(due to Mader \cite{Mader2010}), and it has a negative answer for any $r\geq 8$. We further answer Problem~\ref{seymour} in the negative for $r\in\{6,7\}$. 
Furthermore, we wonder whether every \(r\)-regular digraph \(D\) satisfies that $c(D)=r$ for $r\in\{4,5\}$.

Then we focus on another problem raised by Steiner, which is to determine $\displaystyle \lim_{r\to\infty} \frac{c_r}{r}$. We improve the lower bound from $3/22$ to $1/6$, and the upper bound from $7/8$ to $5/6$. This problem is very interesting and deserves further study.

\subsection*{Acknowledgements}

\begin{sloppypar}

\noindent
This work was supported by the National Key R\&D Program
of China (No. 2023YFA1010202) and the National Natural Science
Foundation of China (Nos. 12471336, 12501473).
\end{sloppypar}

\subsection*{Declaration on the use of generative AI}

\noindent

The counterexamples are constructed entirely by the authors. The authors initially obtained a lower bound on $c_r$
of approximately $0.15r$, and then the authors use AI to optimize several inequalities and improve the lower bound to  $r/6$. The authors are responsible for the mathematical claims and the final text.

\vskip 3mm
\end{spacing}
\end{document}